\documentclass{article}
\usepackage{graphicx}
\usepackage{amsmath}
\usepackage{amssymb}
\usepackage{url}
\usepackage{diagrams}
\usepackage{eepic}
\usepackage{newcent}

\newcommand{\M}{\mathfrak{M}}

\newtheorem{lemma}{Lemma}
\newtheorem{defn}{Definition}
\newtheorem{thm}{Theorem}
\newtheorem{axiom}{Axiom}

\newtheorem{prop}{Proposition}
\DeclareMathOperator{\dist}{dist}
\title{Huygens' principle - a synthetic account}
\author{Anders Kock}

\begin{document}
\date{}

\maketitle

\section*{Introduction}

Huygens' theory of wave fronts in a space $M$ is built on two primitive notions: 1) 
the notion 
of when two submanifolds (of a given manifold $M$) {\em touch} in a 
given point (together with a derived notion of envelope), and 2) the notion of 
{\em sphere} in $M$, with a given length as radius.

The wave front theory deals with the case where the submanifolds are 
hypersurfaces, i.e.\ of dimension $n-1$ (where $M$ is of dimension 
$n$), e.g.\ circles in a plane, or spheres in 3-space.

The principle that carries his name is then: given a hypersurface $B 
\subseteq M$ and a length $s$ (sufficiently small), then each of the spheres $S(b,s)$ with center 
$b\in B$ and radius $s$ touch another hypersurface, the {\em 
envelope} of the spheres $S(b,s)$ as $b$ ranges over $B$. They are  
 ``elementary waves'', or  
``wavelets'' generated by $B$. If we let $B\vdash s$ be the wave 
front obtained this way, it is part if the theory that $(B \vdash 
s)\vdash t = B\vdash (s+t)$; for this part of the theory, one needs 
that one can {\em add} (and to a certain extent, subtract) 
lengths\footnote{
More precisely, we assume that the possible lengths $r$,$s$ etc.\ 
form a commutative cancellative semigroup $R_{+}$, written additively. So if 
$r+s =t$, we say that $r=t-s$; this is well defined, by 
cancellativity. If $r+s=t$, we say that $r<t$ and that $s<t$.}
.
The theory depends on two classical facts:

\medskip

{\em 1) If two (non-concentric) spheres touch, then the point in which they 
touch is unique.

2) Two spheres touch each other (internally) iff the distance between 
their centers is the difference between their radii:}

\bigskip

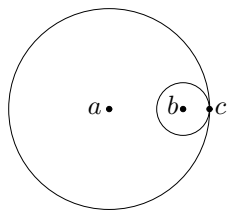
\begin{equation*}\label{dPictx}
\begin{picture}(100,50)(0,0)

\put(30,30){\circle*{2}}
\put(30,30){\circle{76}}
\put(58,30){\circle{20}}
\put(58,30){\circle*{2}}
\put(68,30){\circle*{2}}

\put(22,28){$a$}

\put(52,28){$b$}
\put(70,28){$c$}
\end{picture}
\end{equation*}

\medskip

\noindent (There is a similar fact for {\em external} touching, see 2') below.) Denote the three marked point by $a$, $b$, and $c$, 
respectively:
$a$ and $b$ are the centers of the two spheres in question; and $c$ is 
the point where these spheres touch. In particular, the principle 2) 
is an existence statement: a touching point $c$ {\em exists} (and is 
unique, by 1)). 
 The distance between $a$ and $b$ is 
 the difference $r$ between the radii.
 The middle point $b$ belongs to the  sphere (not depicted) $S(a,r)$ with center $a$ and radius $r$ . In terms of the general picture,
 if the hypersurface $B$ is $S(a,r)$, 
then $S(b,s)$ is one of the wavelets, and $S(a,r+s)$ is the new wave 
front $B\vdash s$, the {\em envelope}  $E$ of the wavelets $S(b,s)$ as $b$ ranges over 
$B$. Both $B$ and some of these wavelets are depicted here:

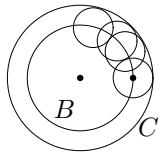
\begin{equation*}\label{cPict3x}
\begin{picture}(100,50)(0,0)
\put(30,30){\circle{40}}
\put(20,15){$B$}
\put(52,8 ){$C$}
\put(30,30){\circle*{2}}
\put(30,30){\circle{55}}
\put(50,30){\circle{15}}

\put(50,30){\circle*{2}}
\put(47,40){\circle{15}}
\put(44,44){\circle{15}}
\put(35,49){\circle{15}}

\end{picture}
\end{equation*}
The envelope $E$ in this case is a sphere $C$. For the general Huygens Principle, $B$ can be any ``hypersurface", at least when $s$ sufficiently small.

There is a fact analogous to 2), for external touching:

\medskip

{\em 2') Two spheres touch each other (externally) iff the distance 
between their centers is the sum of their radii.}

\medskip

An envelope  of a family of subspaces $S_{b}$ of subspaces of a 
space $M$, for  $b\in B$, is a subspace $E$ that touches every 
$S_{b}$, and conversely, every point of $E$ is touched by some 
$S_{b}$. 

Both the notion of touching and the notion of envelope  
are classically formulated in terms of differential calculus. We present an 
alternative way, where these notions
derive from the basic 
primitive relation of (pairs of) {\em neighbour points}, as known in 
synthetic differential geometry (SDG) and in algebraic geometry (where one talks about the {\em first 
neighbourhood of the diagonal} of a scheme $M$, which is a subscheme 
$M_{(1)}\subseteq M\times M$), cf.\ e.g. [SGM] and the references 
therein.

\medskip
Some underlying notion of ``continuity" or ``cohesiveness" of the spaces and of the maps is implicit in all our constructions and arguments; such cohesiveness may be made explicit by saying that everything takes place in a suitable topos, whose objects we call ``spaces", or just ``sets", and which we talk about as if they were sets. The ``elements" of such ``sets", we also sometimes call ``points". This is a standard method in e.g.\ synthetic differential geometry (see e.g.\ \cite{SDG} or \cite{SGM}). For the  specific theory developed in the present article, topos models are described in \cite{MSSDG}. The present text may be seen as a simplification of some of the theory in loc.cit. Also, it differs from it by having the emphasis on ``wave fronts" (hypersurfaces), rather than on ``rays" (geodesics). It may be seen as an attempt to provide a rigourous and constructive account of (part of) the Chapter ``Huygens' Principle" in Arnold's \cite{Arnold} (who states ``we will not pursue rigour here").

\section{Touching}

We are presenting an axiomatic theory for the concepts of touching (and hence for envelopes); 
the novel aspect is the structure that defines the notion of 
{\em touching}. Classically, it  is defined in terms of differential 
calculus; but in the present treatment, touching derives from the basic 
primitive relation of (pairs of) {\em neighbour points}, as known in 
algebraic geometry and SDG. For example, in the category of affine schemes, the {\em first neighbourhood of the 
diagonal} of a scheme $M$ , i.e.\  $M_{(1)}\subseteq M  \times M$ 
defines such a relation, and this is a main model to have in mind.

All the sets we consider are supposed to be equipped with a such a 
neighbour relation $\sim$, and all maps we consider are supposed to 
preserve $\sim$. So what we call ``sets'' we might also be called 
``spaces''.

So a {\em neighbour relation} $\sim$ on a space $M$ is a reflexive 
symmetric relation $M_{(1)}$ on $M$, so $x\sim y$ means $(x,y) \in 
M_{(1)}$.  If $x\sim y$, we say 
that $x$ and $y$ are {\em neighbour points}.

If $C \subseteq M$, the $\sim$ on $M$ restricts to a $\sim$ on $C$.
If $x\in M$, we call the set of points $y$ which are neighbours of 
$x$ {\em the monad} around $x$, and denote it $\M (x)$.

\begin{defn}Let $A$ and $B$ be subspaces of $M$, and let $x\in A \cap 
B$. Then we say that $A$ and $B$ {\em touch each other in $x$}, or that they {\em have $x$ as a touching point}, if for 
all $y\sim x$, we have
$$y\in A \mbox{ iff }y\in B,$$
equivalently if
$$\M(x)\cap A = \M(x)\cap B.$$\end{defn}
If we replace the equality sign here by an inclusion $\subseteq$, we get a weaker 
notion, which is also geometrically relevant, and might be called ``subtouching". 

\medskip

To say that two subsets $A$ and $B$ of $M$ touch in a point $c$ implies 
that $c\in A \cap B$, but  the set of points where $A$ 
and $B$ touch (the ``touching set") is in general smaller than $A\cap B$.

Spheres have the property that the touching set of two spheres is 
a singleton set; this will be one of the axioms below. In the 
analytic model for the axiomatics (cf.\ \cite{MSSDG}), this comes from 
non-singularity of certain second order derivatives.

\section{Spheres}In the following, we consider a fixed ambient space $M$, which one may think of as $n$-dimensional Euclidean space; but in the axiomatic we present, we do not have ``lines", as a primitive notion (rather, the notion of line/geodesic/ray can be defined, essentially following \cite{Busx}, see also \cite{MSSDG}. But the notion of {\em sphere} (special case: {\em circle}) will be crucial:
 
In general terms, a sphere in $M$  is the set of 
points with a given {\em distance} to a given point. 

For 
this to make sense, we need that $M$ is equipped with a {\em metric} 
$\dist$, 
in some suitable sense. For the use of spheres in the Huygens' principle, it suffices with
a quite weak notion of metric (``pre-metric"), which we now make precise. The values 
of the function $\dist$ should be a space $R_{+}$, which one may think of as the strictly 
positive real numbers. For $t$ and $s$ in $R_{+}$, we say that $t>s$ if $t$ is of the 
form $r+s$ for some $r$. We require that this is a total order, in the sense that for any $s$ and $t$ in $R_{+}$, either $r<t$ or $t<r$.
 We further require that $R_{+}$ is a 
{\em cancellative} commutative semigroup. This is here equivalent to: if $t>r$, then there is a {\em unique} $s$ with $r+s=t$. This $s$ may be denoted $t-r$.

By a {\em pre-metric} on $M$, we mean a partially defined 
$R_{+}$-valued function 
$\dist$, with $\dist (a,b)$ defined whenever $a$ and $b$ are {\em distinct} 
points of $M$. There is only one requirement, namely that $\dist$ is 
symmetric, i.e.\ $dist (a,b)= \dist (b,a)$.\footnote{We do not assume any kind 
of triangle inequality, and ``$\dist (a,a)=0$'' does not make sense. 
In fact, if $a\sim b$, then $a$ and $b$ are not distinct, and $\dist 
(a,b)$ is not defined. The reason for considering such weak kind of 
``metric'' is that it is immediately compatible with SDG, see 
\cite{MSSDG}.}

A {\em sphere} $S(a,r)$ in a space $M$ is a subset of $M$ given by its {\em center} 
$a\in M$  and its {\em radius} $r\in R_{+}$, as the set of points in $M$ which 
have distance $r$ to the center $a$,
$$S(a,r):= \{b\in M\mid \dist (a,b)=r\}.$$
Note that the set $S(a,r)$, for a given $a$ and $r$, may well be empty.
But the axioms for touching of spheres in the following section imply existence of touching-points on some spheres.

\section{Touching of spheres}
We give now the axioms which suffice for a proof of the Huygens 
Principle. Two of these axioms have in fact been stated in the Introduction 
as 1), 2) and 2'), 
but now we have made a relevant notion of {\em touching} and {\em sphere} 
explicit. 
\begin{axiom}\label{unix}If two (non-concentric\footnote{meaning: their centers 
are distinct}. For the use of the term ``distinct", see Appendix.) spheres touch, then they have  a 
unique touching point.
\end{axiom}
\begin{axiom}\label{hux}The following two conditions are equivalent, for two 
(non-con\-cen\-tric) spheres:

1) The spheres touch. 

2) The distance between their centers 
is the {\em difference} between their radii, or the distance between their 
centers is the {\em sum} of their radii.
\end{axiom}

The two alternatives in the second condition are mutually exclusive. 
If the first of the alternatives hold, we talk about {\em internal} 
touching, in the second {\em external} touching.

\medskip
 Let $\dist (a,b)$ be denoted $r$, and let $s\in R_{+}$. The point $c$ 
 where $S(a,r+s)$ and $S(b,s)$ touch may be denoted 
 $a\triangleright_{s}b$ (cf.\ \cite{MSSDG}); geometrically, $a\triangleright_{s}b$ is obtained 
 by ``extrapolation'' $s$ units from $b$ along the line given by 
 $a$ and $b$.  
 
 \medskip

From the Axiom (``internal touching''-part) follows that 
every sphere  with center $b$ on $S(a,r)$ and 
radius $s$ touches the sphere $S(a,r+s)$; this is almost saying that 
the sphere $S(a, r+s)$ is an envelope of the family of radius 
$s$-spheres mentioned, which is the simplest case of Huygens' 
priciple. Also, from the Axiom (``external-touching''-part) follows 
that every sphere with radius  $s$ and center on $S(a,r+s)$ touches 
$S(a,r)$. However, an argument is needed in both cases, namely that 
every point of the putative envelope is in fact touched by a sphere 
of the relevant family. It will follow from the 
``reciprocity''-proposition below.

\medskip

For this, and for later use, we need one further, more technical,  axiom (Axiom 1.7 in \cite{MSSDG}):
\begin{axiom}\label{dimx} Assume that $A$ and $C$ are (non-concentric) 
spheres in $M$ 
and that $b\in A\cap C$. Then:  $\M(b) \cap A \subseteq \M(b) \cap 
C$ implies $\M(b) \cap A = \M(b) \cap C$.  
\end{axiom}
This reasonableness of this axiom is related to the fact from linear 
algebra that if 
two linear subspaces of a vector 
space have the same dimension, and one is contained in the other, 
then they are equal. (The Axiom may verbally be rendered: for spheres, ``subtouching implies touching".)

\medskip

Given  three distinct points $a$, $b$,  and $c$, with $\dist (a,b)=r$, $\dist(b,c) =s$, and $\dist(a,c)=r+s$.
\begin{lemma}[Reciprocity] The following conditions  are equivalent. 1) $S(b,s)$ touches $S(a,r+s)$ in $c$; 2) $S(c,s)$ touches $S(a,r)$ in $b$.

\end{lemma} 

In the following picture, the concentric spheres (circles) are $B=S(a,r)$ and 
$C=S(a,r+s)$; the two small spheres (``wavelets'') $S_{1}$ and $S_{2}$ both have radius 
$s$, and have centers $b\in B$ and $c\in C$; $S_{1}$ touches $C$ in 
$c$ (internal touching), $S_{2}$ touches $B$ in $b$ (external touching):

\bigskip

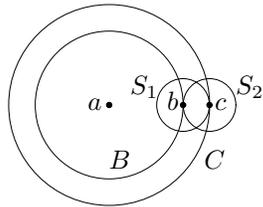
\begin{equation*}\label{dPictx2}
\begin{picture}(100,50)(0,0)

\put(30,30){\circle*{2}}
\put(30,30){\circle{56}}
\put(30,30){\circle{76}}
\put(68,30){\circle{20}}
\put(58,30){\circle{20}}
\put(58,30){\circle*{2}}
\put(68,30){\circle*{2}}

\put(30,6){$B$}
\put(66,6){$C$}
\put(22,28){$a$}
\put(38,34){$S_{1}$}
\put(78,34){$S_{2}$}

\put(52,28){$b$}
\put(70,28){$c$}
\end{picture}
\end{equation*}

\medskip

\noindent {\bf Proof.} With the notation from the figure, there is a map $B\to C$ which to any $b\in B$ associates the unique point where $S(b,s)$, according to the Axioms 1 and 2, touches $C$; geometrically, this map is ``radial projection" (in the notation referred to above, this is the map  $b\mapsto a\triangleright_{s}b$). Likewise, there is a map $C\to B$ which to any $c\in C$ associates the unique point where 
$S(c,s)$ touches $B$.
The assertion of the Lemma is that these two maps are inverses of each other. 
We prove that the composite $ B\to C \to B$ is the identity (the proof that the composite $C\to B\to C$ 
is the identity map is similar). So for $b\in B$, let $c$ be the point where $S(b,s)$ touches $C$; 
to prove that $S(c,s)$ touches $B$ in $b$, it suffices by Axiom 3 to prove that $\M (b)\cap R \subseteq S(c,s)$. 
If $b'$ in $R$  has $b'\sim b$, then its image $c'$ under the central projection map 
$B\to C$ is $\sim c$, and by construction of $c'$ as a touching point on $S(b',s)$, 
we have not only $\dist (b',c')=s$, but $\dist (b',c'')=s$ for all $c''\sim c'$ in $C$. 
This applies in particular if we take  $c''=c$; so we conclude that $\dist (b',c)=s$; so $b'\in S(c,s)$, as desired.
(Note that we only used the touching point property of $c$ to make 
sure that $c'\sim c$; for, then they are the respective images of a 
certain map $B\to C$.)

\medskip

\medskip



\section{Contact elements}

A {\em contact element} $P$ in $M$ at\footnote{the phrase ``at $b$'' 
is redundant under the slightly stronger axiomatics of \cite{MSSDG}, 
because there, the $x$ can be reconstructed from the set $P$, namely 
as its focus, in the sense of loc.cit.} $b\in M$ is a set of the form $\M(b)\cap B$, 
where $B$ is a sphere containing $b$; we then say that the sphere $B$ 
{\em touches $P$ in} $b$, or just that it touches $P$, when $b$ is understood. Two spheres containing $b$ 
define the same contact element (in $b$) if they touch each other (in $b$). 

Let $P$ be a contact element at $b$, and let $x$ be distinct from $b$, 
hence also distinct from any $b'\sim b$, in particular, distinct 
from all points of $P$. We say that $x\perp P$ (read: $x$ is {\em orthogonal to } $P$) if $\dist (b',x)=\dist 
(b,x)$ for all $b'\in P$. This is to say that $\M(b)\cap B \subseteq S(x,s)$, where $s=\dist(b,x)$; by Axiom 3, this is again equivalent to saying that the spheres $B$ and $S(x,s)$ touch in $b$, and again equivalent to: the sphere $S(x,s)$ 
touches $P$.

The spheres touching $P$ fall in two classes, where two spheres $B_{1}$ and $B_{2}$ 
are in the same class if they touch each other from the inside. 
Selecting one of these classes gives what may be called a 
(transversal) {\em orientation} of $P$; the spheres in the selected 
class are then said to touch $P$ {\em internally} or from the {\em inside}, and those in the 
other class are said to touch $P$ {\em externally}  or from the {\em outside}.

 If $P$ is given an orientation, we say that an $x$ (with $x\perp P$) is on the {\em 
positive} side of $P$ if $S(x,s)$ touches $P$ on the outside (where $s=\dist(b,x)$). 

\begin{prop}Let $P$ be an oriented contact element at $b$, and let $s\in R_{+}$. Then there is exactly one point $c$ such that $c\perp P$, with $c$ on the positive side of $P$ and with  $\dist(b,c)=s$.
\end{prop}
{\bf Proof.} Pick a sphere touching $P$ from the inside, say $P=\M(b)\cap S(a,r)$. Then $c$ is  outside of $S(a,r)$. The assumed property of $c$ states  that
\begin{equation}\label{cicx}\M(b)\cap S(a,r) \subseteq S(c,s),\end{equation} which in view of Axiom 3 implies that $S(a,r)$ and $S(c,s)$ touch in $b$.
By Axiom 2, this means that either $\dist(a,c)$ is either $ r+s$ or  $r-s$. 
The last case is excluded since $c$ is outside the sphere $S(a,r)$. Applying the Reciprocity Lemma then gives that $c$ is the touching point of  
$S(a,r+s)$ and $S(b,s)$, and the uniqueness of $c$ then follows since the assumption on $c$ only refers to the oriented contact element, not to its presentation by a particular sphere like $S(a,r)$.

\medskip

The unique point  $c$ described in this Proposition deserves a notation, we write $P\vdash s$ for it (more fully $(P,b)\vdash s$); it is thus characterized  by: $c\perp P$, $\dist(b,c)=s$, and being on the positive side of $P$. The proof gives that it may be {\em constructed} by  {\em picking} an arbitrary sphere $S(a,r)$ touching $P$ from the inside. In terms of such sphere
$c$ is characterized by being outside $S(a,r)$ and satisfying
\begin{equation}\label{condx}\forall b' \;  [(b'\sim b \wedge \dist(a,b')=r )\Rightarrow (\dist(b',c)=s)].\end{equation}

\medskip

Note that the condition that $S(a,r)$ touches $S(c,s)$ in $b$ implies that (\ref{cicx}) holds, equivalently
$$\forall b': [b'\in P \Rightarrow \dist (b',c)=s]$$
which may be written 
\begin{equation}\label{charptx}c\in \bigcap _{b'\in P}S(b',s);\end{equation}
The condition (\ref{charptx}) does not involve any somewhat arbitrarily chosen sphere $S(a,r)$. It shows that the point $P\vdash s$ is a {\em characteristic point}, in the discriminant-sense (cf.\ \cite{END} for a discussion of this ``predicative" way of constructing envelopes).


\section{Hypersurfaces, and Huygens' Principle}
By a {\em hypersurface} in $M$, we understand a subset $B$ such that for every $b\in B$, the set $\M (b)\cap B$ 
is a contact element at $b$, which we then denote $B(b)$. To give $B$ an  {\em orientation} is to give 
each of these contact elements an orientation. Let $s\in R_{+}$. Then we have 
the contact element $B(b)$, 
and the point $B(b)\vdash s$. We have a map $B \to M$ sending $b\in B$ to $B(b)\vdash s$. For $s$ suitably small, we assume that this map is injective.  (If for instance $B$ is the sphere $S(a,r)$, then either $s$ can be chosen freely, or $s$ has to be $<r$, depending on the orientation given $S(a,r)$. For the orientation where the concave side is the positive one, $B\vdash r$ equals $a$ for all $b\in B$, so $B\vdash -$ cannot be injective.)

We consider in the following a fixed oriented hypersurface $B$ and an $s$ small enough so that the map $B\to M$ thus described is injective; its image is conveniently denoted $B\vdash s$. It is the set $C$ 
of points of the form $B(b) \vdash s$ for $b\in B $. By assumption,
the map $b\mapsto B(b)\vdash s$ is a bijection $B\to C$. The point $b\in B$ corresponding to $c\in C$ under this bijection is what classically is called the {\em foot} of $c$ on $B$. 
We need to assume a certain ``continuity" property of the situation: if a point $c$ is of the form $B(b)\vdash s$, then any point $x\sim  c$ is likewise of the form $B(b')\vdash s$ for some $b'\sim b$ and some $s'\in R_{+}$. In particular, $x$ is on the positive side of $B(b')$.

Under these circumstances:

\begin{thm}[Huygens' principle] The set $B\vdash s$ is a hypersurface. It is an envelope of the family 
$\{S(b,s)\mid b\in B \}$.
\end{thm} 
{\bf Proof.} 
Consider a $c\in C$, say $c=B(b)\vdash s$ (so $b$ is the foot of $c$ on $B$).  We shall prove that
\begin{equation}\label{mmx}
\M(c)\cap C = \M(c)\cap S(b,s).
\end{equation}
To prove the inclusion $\M(c)\cap C \subseteq \M(c)\cap S(b,s)$, let $c'\in C$ with $c'\sim c$. 
So $c'$ is of the form $B(b')\vdash s$ for a unique $b'\in B$ (the foot of $c'$ on $B$). So $c'\perp B(b')$ with $\dist (b',c')=s$. Since $b\in B(b')$, 
we therefore have $\dist(b,c')=s$, so $c'\in S(b,s)$.

To prove the inclusion $\M(c)\cap S(b,s)\subseteq \M(c)\cap C $, let $x\sim c$ and $\dist(b,x)=s$.
Let $b'$ be the foot of $x$ on $B$, i.e.\ $x=B(b')\vdash s'$ .Then $x'\sim c$ implies $b'\sim b$. 
Since $x\perp B(b')$ and $b\in B(b')$, we have 
$\dist (b',x) = \dist (b,x) =s$. But $x\perp B(b')$ and $\dist(b',x)=s$ characterize $B(b')\vdash s$, 
so $x=B(b')\vdash s$, and thus $x\in C$.
 
These two inclusions prove 
(\ref{mmx}), 
and since the right hand side of it is a contact element, then so is the left hand side. So $C$ is a hypersurface, proving the first assertion.
 For the second assertion, given $c\in C$. By construction of $C$, $c$ is of the form $B(b)\vdash s$. Then (\ref{mmx}) proves that $S(b,s)$ touches $C$ in $c$.  Conversely, given $b\in B$, take $c:=B(b)\vdash s$; then the same equality shows that every sphere of the family of $S(b,s)$'s touches $C$. The correspondence between $b\in B$ and $c\in C$ is bijective.
 
 \section*{Appendix on logic} As in all situations where continuity, or some other cohesiveness, is presupposed (built into the underlying category), the law of excluded middle has limited applicability. So even when $a$ is {\em not} distinct (distinguishable) from $b$, one cannot conclude that $a=b$.
 
 In the present case, one should note the particular relationship between the {\em neighbour} notion and the notion of {\em distinct} points: {\em Neighbour points are never distinct.}) If $a$ is distinct from $b$, then also any neighbour of $a'$ of $a$ is distinct from $b$.
 
 The models (cf.\ \cite{MSSDG}) of the present axiomatic theory usually presuppose a number line $R$, which has the structure of a commutative ring. In such a ring, $a$ is {\em distinct} from $b$ if $b-a$ is multiplicatively invertible. And $b'	$ is a (first order) {\em neighbour} of $b$ if $(b'-b)^2=0$. (One also has a notion of higher order neighbours: $(b'-b)^n =0)$.

\medskip

\noindent Anders Kock
 
 \noindent
 Dept.\ of Mathematics, \\ University of Aarhus, Denmark
 
 \noindent kock (at) math. au.dk
 
 \noindent April 2018

 \end{document}